\documentclass[10pt]{amsart}

\usepackage[latin1]{inputenc}
\usepackage{amsmath}
\usepackage{amsfonts}
\usepackage{amssymb}
\usepackage[mathscr]{eucal}
\usepackage{diagrams}
\usepackage{xy,color,graphicx}
\usepackage{hyperref}
\xyoption{all}

\newcommand{\wis}[1]{{\text{\em \usefont{OT1}{cmtt}{m}{n} #1}}}

\newcommand{\C}{\mathbb{C}}

\newcommand{\mon}{\mathbb{N}^{\times}_+}

\newtheorem{definition}{Definition}

\newtheorem{theorem}{Theorem}
\newtheorem{corollary}{Corollary}

\newtheorem{remark}{Remark}

\title{Supernatural Numbers and a new Topology on the Arithmetic Site}
\author{Lieven Le Bruyn} 
\address{Department of Mathematics, University of Antwerp \\ 
 Middelheimlaan 1, B-2020 Antwerp (Belgium) \\ {\tt lieven.lebruyn@uantwerpen.be}}

\begin{document}
\sloppy

\maketitle

Alain Connes and Caterina Consani introduced in \cite{CC} and \cite{Connes} the {\em arithmetic site} as the topos $\widehat{\mon}$ of sheaves of sets on the small category corresponding to the multiplicative monoid $\mon$, equipped with the chaotic topology. They proved that the isomorphism classes of points of this topos are in canonical bijection with the finite ad\`ele classes
\[
\wis{points}(\widehat{\mon}) = \mathbb{Q}^{\ast}_+ \backslash \mathbb{A}^f_{\mathbb{Q}} / \hat{\mathbb{Z}}^{\ast} \]
The induced topology on this set, coming  from the locally compact topology on the finite ad\`ele
 ring $\mathbb{A}^f_{\mathbb{Q}}$, is trivial, and therefore this 'space' is best studied by the tools of noncommutative geometry, as was achieved in \cite{CC}. 
 
However, we will define another topology on this set of points. Even though the {\em SGA4-topology}, see\cite[IV.8]{SGA4}, on the arithmetic site is also trivial, we can define for any sieve $C$ in $\mon$ 
\[
\mathbb{X}(C) = \wis{points}(\widehat{\mon}) \cap \wis{points}(\widehat{C}) \]
The corresponding {\em sieve-topology} on the points $\mathbb{Q}^{\ast}_+ \backslash \mathbb{A}^f_{\mathbb{Q}} / \hat{\mathbb{Z}}^*$ shares many properties one might expect the mythical object $\overline{\wis{Spec}(\mathbb{Z})}/\mathbb{F}_1$ to have: it is compact, has an uncountable basis of opens, each nonempty open being dense, and, it satisfies the $T_1$ separation property for incomparable points.

\section{Arithmetic sites}

In this section we translate topos-theoretic terminology to the setting of interest. With $\mon$ we denote the multiplicative monoid of all strictly positive natural numbers. We define a class of sub-monoids of $\mon$, which will correspond to the sieves of the associated small category, whence the s-terminology.

\begin{definition} A multiplicative sub-monoid $\wis{C}$ of $\mon$ is said to be an
\begin{enumerate}
\item{{\em s-monoid} if for all $c \in \wis{C}$ also all additive multiples $c.n \in \wis{C}$.  Every s-monoid has a presentation
\[
\wis{C}=\wis{C}(c_1,c_2,\hdots) = c_1 \mon \cup c_2 \mon \cup \hdots \]
for possibly infinitely many elements $c_i$ from $\wis{C}$.}
\item{{\em S-monoid} if it is an s-monoid which is also additively closed, that is, for all $c,c' \in \wis{C}$ we have $c+c' \in \wis{C}$. Every S-monoid has a presentation
\[
\wis{C} = \wis{C}_+(c_1,c_2,\hdots,c_k) = (c_1.\mathbb{N}+c_2.\mathbb{N}+\hdots+c_k.\mathbb{N}) \cap \mon\]
for finitely many elements $c_i$ from $\wis{C}$.}
\end{enumerate}
\end{definition}

To an s-monoid $\wis{C}$ we associate a small category, also denoted $\wis{C}$, which has just one object $\bullet$ with its monoid of endomorphisms, $\wis{C}(\bullet,\bullet)$, under composition isomorphic, as monoid, to the multiplicative sub-monoid $\wis{C} \cup \{ 1 \}$ of $\mon$ (here, $1$ corresponds to the obligatory identity morphism $id_{\bullet}$). If, for $\wis{C} \not= \mon$, this small distinction may cause confusion, we will use the adjectives 'monoid' or 'category' to distinguish between the two uses of $\wis{C}$.

\vskip 2mm

From \cite[p. 37]{MM} we recall that a {\em sieve} in the category $\wis{C}$ is a collection $\wis{S}$ of arrows in $\wis{C}$ with the property that if $f \in \wis{S}$ and if the composition $f \circ h$ exists in $\wis{C}$ then also $f \circ h \in \wis{S}$. That is, a sieve $S$ is of the form
\[
\wis{S} = x_1.\wis{C} \cup x_2.\wis{C} \cup \hdots = \wis{S}(x_1,x_2,\hdots) \]
The set $\Omega_{\wis{C}}$ of all sieves on $\wis{C}$ is a lattice (under $\cap$ and $\cup$) and is ordered (under $\subset$) having a unique maximal element $\wis{C}$.

Observe that the arrows in a sieve $\wis{S}$ of $\mon$ is the same thing as an s-monoid in $\mon$. Moreover, the same holds for any sieve $\wis{S}$ of the category $\wis{C}$ corresponding to an s-monoid $\wis{C}$.

\vskip 2mm

The set $\Omega_{\wis{C}}$ has a right action by $\wis{C} \cup \{ 1 \}$ by $S \odot c = c^{-1}.S \cap (\wis{C} \cup \{ 1 \})$.
A {\em Grothendieck topology}, see \cite[III.2]{MM}, on the category $\wis{C}$ is a subset $\mathcal{J} \subset \Omega_{\wis{C}}$ satisfying the following properties
\begin{enumerate}
\item{$\wis{C} \in \mathcal{J}$}
\item{(stability) if $\wis{S} \in \mathcal{J}$, then $\wis{S} \odot \wis{C} = \{ \wis{S} \odot c~|~c \in \wis{C} \} \subset \mathcal{J}$}
\item{(transitivity) if $\wis{S} \in \mathcal{J}$ and if $\wis{R} \in \Omega_{\wis{C}}$ such that $\wis{R} \odot \wis{S} = \{ \wis{R} \odot s~|~s \in \wis{S} \} \subset \mathcal{J}$, then $\wis{R} \in \mathcal{J}$.}
\end{enumerate}
Observe that it follows that if $\wis{S} \in \mathcal{J}$ and $\wis{S} \subset \wis{S'}$ in $\mathcal{S}$, then also $\wis{S'} \in \mathcal{J}$.  The coarsest Grothendieck topology on $\wis{C}$, with $\mathcal{J}_{ch} = \{ \wis{C} \}$, is called the {\em chaotic topology} on the category $\wis{C}$, see \cite[II.1.1.4]{SGA4}.

\vskip 2mm

A {\em presheaf} on the category $\wis{C}$, equipped with a Grothendieck topology $\mathcal{J}$, is a contravariant functor
\[
P~:~\wis{C} \rTo \wis{Sets} \]
and hence corresponds to a set $R=P(\bullet)$, equipped with a {\em right}-action of the monoid $\wis{C}$. We will then denote the pre-sheaf by $P_R$, and we have an equivalence of categories between the category of all presheaves $\wis{PreSh}(\wis{C},\mathcal{J})$, with natural transformations as morphisms, and, the category $\wis{Sets}-\wis{C}$ of all right $\wis{C}$-sets with $\wis{C}$-maps as morphisms.

\vskip 2mm

The {\em Yoneda embedding} $\wis{y}~:~\wis{C} \rTo \wis{PreSh}(\wis{C},\mathcal{J})$ sends $\bullet$ to the functor $y(\bullet)=\wis{C}(-,\bullet)$, that is, to the pre-sheaf $P_{\wis{C}}$ where the monoid $\wis{C}$ has the right-action by multiplication on itself. As such, sieves can be viewed as sub-functors of $\wis{y}(\bullet)$.

\vskip 2mm

A presheaf $P=P_R \in \wis{PreSh}(\wis{C},\mathcal{J})$ is said to be a {\em sheaf} if and only if for every sieve $\wis{S} \in \mathcal{J}$, the sub-functor $\wis{S} \rInto \wis{y}(\bullet)$ induces an isomorphism between the natural transformations $\wis{Nat}(\wis{S},P_R) \simeq \wis{Nat}(\wis{y}(\bullet),P_R)$, see \cite[III.4]{MM}. 

In terms of the right $\wis{C}$-set the sheaf property says that there is a natural isomorphism
\[
\wis{Maps}_{\wis{C}}(\wis{S},R) = \wis{Maps}_{\wis{C}}(\wis{C},R) \]
That is, every right $\wis{C}$-map $\wis{S} \rTo R$ extends uniquely to a right $\wis{C}$-map $\wis{C} \rTo R$. 

This can be reformulated in terms of {\em matching families}, see \cite[III.4]{MM} as follows: if we have a family of elements $r_x \in R$ for all $x \in \wis{S}$ such that $r_{x.c} = r_x.c$ for all $c \in \wis{C}$, then there is a unique element $r \in R$ such that each $r_x = r.x$. 

The category of all sheaves on $\wis{C}$ for the Grothendieck topology $\mathcal{J}$ will be denoted by $\wis{Sh}(\wis{C},\mathcal{J})$. Categories equivalent to sheaf categories are called {\em topoi}.

If we take the chaotic topology $\mathcal{J}_{ch} = \{ \wis{C} \}$, then the sheaf-condition is void whence $\wis{PreSh}(\wis{C},\mathcal{J}_{ch} ) \simeq \wis{Sh}(\wis{C},\mathcal{J}_{ch}) \simeq \wis{Sets}-\wis{C}$.

\begin{definition} The {\em arithmetic site} corresponding to the s-monoid $\wis{C}$ is the topos of all sheaves on the category $\wis{C}$, equipped with the chaotic topology, and will be denoted by
\[
\widehat{\wis{C}} \simeq \wis{Sh}(\wis{C},\mathcal{J}_{ch}) \simeq  \wis{Sets}-\wis{C} \]
\end{definition}

\section{The points}

In this section we will characterize the isomorphism classes of points in the arithmetic sites $\widehat{\wis{C}}$. We begin by motivating the definition. 

If $X$ and $Y$ are (ordinary) topological spaces and if $X \rTo^f Y$ is a continuous map, then there is an adjoint pair between the sheaf-topi
\[
\xymatrix{\wis{Sh}(X) \ar@/^2ex/[rr]^{f_*} & & \wis{Sh}(Y) \ar@/^2ex/[ll]_{f^*}} \]
Here $f_*$ is the "direct image" functor, $f^*$ is the "inverse image" functor, and $f^*$ is a left adjoint to $f_*$. Further, $f^*$ preserves finite limits, thus $f^*$ is left exact, see \cite[p. 348]{MM}. In particular, if $X$ is a one-point topological space $\{ p \}$, then $\wis{Sh}(X) \simeq \wis{Sets}$ and so the embedding $f$ gives rise to a {\em geometric morphism}, that is an adjoint pair of functors between $\wis{Sets}$ and $\wis{Sh}(Y)$ which motivates the definition of points of a topos, see \cite[VII.5]{MM}.

\begin{definition} A {\em point} $p$ of the arithmetic site $\widehat{\wis{C}}$ is a geometric morphism, that is, a pair of functors
\[
\xymatrix{\wis{Sets} \ar@/^2ex/[rr]^{f_*} & & \widehat{\wis{C}} \ar@/^2ex/[ll]_{f^*}} \]
such that $f^*$ is a left adjoint to $f_*$ and $f^*$ is left-exact. Here, the functor $f^*$ can be viewed as taking the stalk of a sheaf at $p$, and, the functor $f_*$ assigns to a set, the corresponding  'skyscraper sheaf' at $p$.
\end{definition}

The functor $f^*$ is determined up to isomorphism by its composition with the Yoneda embedding $\wis{y}$, that is, by the functor
\[
\wis{C} \rTo^{\wis{y}} \widehat{\wis{C}} \rTo^{f^*} \wis{Sets} \]
and hence we have to consider suitable {\em covariant} functors $A~:~\wis{C} \rTo \wis{Sets}$. Such functors correspond to sets $L = A(\bullet)$, this time with a {\em left} $\wis{C}$-action and we will denote the functor $A_L$. The functor $f^* : \widehat{\wis{N}} \rTo \wis{Sets}$ corresponding to $A_L$ is given by
\[
\widehat{\wis{C}} \rTo^{f^*} \wis{Sets} \qquad P_R \mapsto R \otimes L \]
where $R \otimes L$ is the quotient-set of the product $R \times L$ modulo the equivalence relation induced by all relations $(r.c,l) = (r,c.l)$ for all $r \in R, l \in L$ and all $c \in \wis{C}$, see \cite[p. 379]{MM}. Such functors have an adjoint functor $f_* : \wis{Sets} \rTo \widehat{\wis{C}}$ determined by sending a set $S$ to the $\wis{C}$-set
\[
S \mapsto \wis{Maps}(L,S) \]
where the right action is given by $\phi.c = \phi \circ (c.-)$. This leaves us to determines the left $\wis{C}$-sets $L$ which are {\em flat}, that is, such that $f^* = - \otimes L$ is left exact.

\begin{theorem} A point $p$ of the arithmetic site $\widehat{\wis{C}}$ corresponds to a non-empty left $\wis{C}$-set $L$ satisfying the following properties:
\begin{enumerate}
\item{$\wis{C}$ acts freely on $L$, that is, for all $a \in L$ and all $c,c' \in \wis{C}$, if $c.a=c'.a$ then $c=c'$.}
\item{$L$ is of rank one, that is, for all $a,a' \in L$ there is an element $b \in L$ and elements $c,c' \in \wis{C}$ such that $a=c.b$ and $a'=c'.b$.}
\end{enumerate}
\end{theorem}

\begin{proof}
This follows by applying Grothendieck's construction of filtering functors to characterize flat functors, see \cite[VII.6, Thm. 3]{MM}. 
\end{proof}

We are interested in the isomorphism classes of points of $\widehat{\wis{C}}$, that is, in left $\wis{C}$-sets $L$ satisfying the requirements of the theorem, upto isomorphism as left $\wis{C}$-set. The next result shows that we can always realize every point upto isomorphism as a specific subset of the strictly positive rational numbers $\mathbb{Q}_+$:

\begin{theorem} Every point $p$ of the arithmetic site $\widehat{\wis{C}}$ for an s-monoid $\wis{C}$ is isomorphic to a subset $L \subset \mathbb{Q}_+$ of the form
\[
L = L(c_1,c_2,\hdots) = \wis{C} \cup \bigcup_{i=1}^{\infty} \wis{C}.\frac{1}{c_i} \]
for certain elements $c_i \in \wis{C}$ satisfying the divisibility condition $c_1 | c_2 | c_3 | \hdots$. 
Conversely, any such subset is indeed a point of $\widehat{\wis{C}}$.
\end{theorem}

\begin{proof} Let $p$ correspond to the left $\wis{C}$-set $L$ satisfying the properties of the previous theorem. Take $l_0 \in L$ and send $l_0 \mapsto 1$, then we have an isomorphism as $\wis{C}$-sets between $\wis{C}.l_0$ and $\wis{C} \subset \mathbb{Q}_+$. If $l'_0 \in L - \wis{C}.l_0$, then there exists an element $l_1 \in L$ and elements $c_1,c'_1 \in \wis{C}$ such that $l_0 = c_1.l_1$ and $l'_0 = c'_1.l_1$ and send $l_1 \mapsto \frac{1}{c_1}$, then the left $\wis{C}$-subset of $L$ spanned by $l_1$ (which contains $l_0,l_0'$) is isomorphic as $\wis{C}$-set to $C.\frac{1}{c_1}$.

If there is an element $l'_1 \in L-\wis{C}.l_1$, take an element $l_2$ and elements $d_2,d'_2 \in \wis{C}$ such that $l_1 = d_2.l_2$ and $l'_1 = d'_2.l_2$ and send $l_2 \mapsto \frac{1}{c_2}$ where $c_2=c_1.d_2$ and the $\wis{C}$-subset of $L$ spanned by $l_2$ (which contains $l_0,l'_0,l_1,l'_1$) is isomorphic to $\wis{C}.\frac{1}{c_2}$. Iterating this process we obtain the claimed isomorphism. 

Conversely, $\wis{C}$ acts clearly freely on any such subset, and if $\frac{a}{c_i}$ and $\frac{b}{c_j}$ are arbitrary elements of this set (with $a,b \in \wis{C}$) and if $i \leq j$, then we can consider the element $\frac{1}{c_j}$ and have elements $b$ and $b' = a.\frac{c_j}{c_i}$ (both in $\wis{C}$ because it is an s-monoid)
to satisfy the rank one condition.
\end{proof}

\begin{corollary} If $\wis{C}$ is an S-monoid, every point in $\widehat{\wis{C}}$ is isomorphic to the positive part of an additive subgroup of $\mathbb{Q}$.
\end{corollary}

\begin{proof}
If $\wis{C}$ is an S-monoid, it is additively closed, and hence, so are each of the subsets $\wis{C}.\frac{1}{c_i}$ of $\mathbb{Q}_+$, proving the claim.
\end{proof}

A {\em supernatural number} (also called a {\em Steinitz number}) is a formal product $s = \prod_{p \in \mathbb{P}} p^{s_p}$ where $p$ runs over all prime numbers $\mathbb{P}$ and each $s_p \in \mathbb{N} \cup \{ \infty \}$. The set $\mathbb{S}$ of all supernatural numbers forms a multiplicative semigroup with multiplication defined by exponent addition and the multiplicative semigroup $\mon$ embeds in $\mathbb{S}$ via unique factorization.

If $s \in \mathbb{S}$, we will write $s_p$ for the exponent of $p \in \mathbb{P}$ appearing in $s$. The {\em support} of $s$ is the set of prime numbers $p$ for which $s_p > 0$. For two supernatural numbers $s$ and $s'$ we say that $s$ divides $s'$, written $s \mid s'$ if $s'=s.s"$ for some supernatural number $s"$, or equivalently, if for all primes $p$ we have $s_p \leq s_p'$.

\vskip 2mm

In \cite[Thm.1,Thm.2]{BZ} Beaumont and Zuckerman characterize all additive subgroups of $\mathbb{Q}$ as being determined by a couple
\[
(a,s) \in \mathbb{N}_+ \times \mathbb{S} \]
with $a$ being coprime to all primes $p$ in the support of $s$. The corresponding additive subgroup is then
\[
\mathbb{Q}(a,s) = \{ \frac{a.n}{m}~|~n \in \mathbb{Z}, m \in \mathbb{N}_+,~m \mid s \} \]
Clearly, $\wis{C}$ acts freely on the positive part $\mathbb{Q}_+(a,s)$ of $\mathbb{Q}(a,s)$, and, two of these are isomorphic, as $\wis{C}$-set, if and only if the corresponding additive groups are isomorphic, that is, if there is a positive rational number $q \in \mathbb{Q}^*_+$ such that $q.\mathbb{Q}_+(a,s) = \mathbb{Q}_+(a',s')$.

That is, up to isomorphism, we need only consider the positive parts $\mathbb{Q}_+(s)= \mathbb{Q}_+(1,s)$, and we have $\mathbb{Q}_+(s) \simeq \mathbb{Q}_+(s')$ if and only if $s$ and $s'$ are equivalent under
\[
s \sim s'~\text{iff}~\begin{cases} s_p = \infty \Leftrightarrow s'_p = \infty \\
s_p = s'_p~\text{for all but at most finitely many $p$} \end{cases} \]

\begin{remark} This equivalence relation on the supernatural numbers $\mathbb{S}$ also appears in the noncommutative geometry of $C^*$-algebras. Recall that a uniformly hyperfinite, or UHF, algebra $A$ is a $C^*$-algebra that can be written as the closure, in the norm topology, of an increasing union of finite-dimensional full matrix algebras
\[
M_{c_1}(\C) \rInto M_{c_2}(\C) \rInto M_{c_3}(\C) \rInto \hdots = A \]
By the double centralizer result it follows that $c_1 | c_2 | c_3 \hdots$. James Glimm proved in \cite{Glimm} that the supernatural number $s=\prod_i c_i$ is an isomorphism invariant for $A$ among UHF-algebras. The algebra $A_{\infty}$ corresponding to $s_{\infty} = \prod_{p \in \mathbb{P}} p^{\infty}$ is often called the universal UHF-algebra.

A {\em noncommutative space} is a Morita equivalence class of $C^*$-algebras. In \cite{Effros} it is shown that two UHF-algebras $A$ and $A'$ with corresponding invariants $s,s' \in \mathbb{S}$ are Morita equivalent if and only if $s \sim s'$, and that the additive subgroup $\mathbb{Q}(s)$ is the Grothendieck group $K_0(A)$, so its positive part $\mathbb{Q}_+(s)$ can be seen as the positive cone of $K_0(A)$. Hence, we can view the set of equivalence classes
\[
\mathbb{S}/\sim~\simeq \wis{moduli}(\wis{UHF}) \]
as the moduli space for noncommutative spaces corresponding to UHF-algebras.

All strictly positive integers $n$ lie in the equivalence class of $1$, expressing the fact that $M_n(\C)$ is Morita equivalent to $\C$, or that the Brauer group $Br(\C)$ is trivial. In this way we can view $\mathbb{S}/\sim$ with the induced well defined multiplication as the Brauer-monoid $Br_{\infty}(\C)$ of UHF-algebras. 
\end{remark}

The set of equivalence classes $\mathbb{S}/\sim$ is also in canonical bijection with the set of finite ad\`ele classes $\mathbb{Q}^{\ast}_+ \backslash \mathbb{A}^f_{\mathbb{Q}} / \hat{\mathbb{Z}}^{\ast}$.
Indeed, for the profinite integers $\widehat{\mathbb{Z}} = \prod_{p \in \mathbb{P}} \widehat{\mathbb{Z}}_p$ and the finite ad\`ele ring $\mathbb{A}^f_{\mathbb{Q}} = \mathbb{Q} \otimes \widehat{\mathbb{Z}}$ we have a canonical bijection between $\mathbb{S} \leftrightarrow \widehat{\mathbb{Z}}/\widehat{\mathbb{Z}}^*$ and $\mathbb{A}^f_{\mathbb{Q}} / \widehat{\mathbb{Z}}^* \leftrightarrow \mathbb{S}_f$ where $\mathbb{S}_f$ are the fractional supernatural numbers, that is, formal products $f=\prod p^{f_p}$ where we allow finitely many $f_p$ to be strictly negative integers. If we take $a = \prod_{p:f_p<o}p^{-f_p}$ then $f$ corresponds to $(a,s) \in \mathbb{N}_+ \times \mathbb{S}$ with $s=a.f$ and we have that $a$ is coprime with the support of $s$. Thus, $\mathbb{A}^f_{\mathbb{Q}} / \widehat{\mathbb{Z}}$ classifies additive subgroups of $\mathbb{Q}^+$ and the $\mathbb{Q}^*$-action corresponds to taking isomorphism classes of such subgroups,which we know are the equivalence classes of supernatural numbers.

\vskip 2mm

We will now characterize which positive cones $\mathbb{Q}_+(s)$ are points of the topos $\widehat{\wis{C}}$ when $\wis{C}$ is an s-monoid. If $n \in \mon$, then $n \sim 1$ and as $\mathbb{Q}_+(1) = \mon$, the only s-monoid for which $\mathbb{Q}_+(1)$ satisfies the rank one condition is $\mon$. Hence,  the only topos for which $\mathbb{Q}_+(n)$ is a point is $\widehat{\mon}$. Remains to deal with the $s \in \mathbb{S}-\mon$:

\begin{theorem} For $\wis{C}$ an s-monoid properly contained in $\mon$, the positive cone $\mathbb{Q}_+(s)$, for $s \in \mathbb{S}-\mon$. is a point of the topos $\widehat{\wis{C}}$  if and only if there exist elements $c_1,c_2,\hdots \in \wis{C}$ such that the supernatural number $\prod_{i=1}^{\infty} c_i$ divides $s$.
\end{theorem}

\begin{proof} Assume $\prod_{i=1}^{\infty} c_i | s$, we have to verify that $\mathbb{Q}_+(s)$ satisfies the rank one condition with respect to $\wis{C}$. Again, because $\wis{C}$ is an s-monoid, it suffices to verify this for all additive generators $\frac{1}{n}$ and $\frac{1}{m}$ of $\mathbb{Q}_+(s)$ with natural numbers $n$ and $m$ dividing $s$. But then, also $lcm(m,n) | s$. As $\prod_{i=1}^{\infty} c_i | s$ there must be an element $c_i \in \wis{C}$ such that also $N=c_i.lcm(m,n) | s$. But then $\frac{1}{N} \in \mathbb{Q}_+(s)$ and we have, if $lcm(m,n)=a.n=b.m$ that
\[
\frac{c_i.a}{N} = \frac{1}{n} \quad \text{and} \quad \frac{c_i.b}{N} = \frac{1}{m} \]
Conversely, assume that $\mathbb{Q}_+(s)$ is a point of $\widehat{\wis{C}}$ and let $n$ be a natural number dividing $s$. Let $a \in \mon - \wis{C}$, then because $\frac{1}{n}$ and $\frac{a}{n}$ belong to $\mathbb{Q}_+(s)$ there must be a natural number $N | s$ and elements $c_1,c'_1 \in \wis{C}$ such that
\[
\frac{c_1}{N} = \frac{1}{n} \quad \text{and} \quad \frac{c'_1}{N} = \frac{a}{n} \]
but then $c_1.n=N$ and so $c_1 | s$. Start again with the elements $\frac{1}{c_1}$ and $\frac{a}{c_1}$ in $\mathbb{Q}_+(s)$ to get $M | s$ and elements $c_2,c'_2$ giving $\frac{c_2}{M} = \frac{1}{c_1}$, to obtain $c_1.c_2 = M | s$. Iterating this procedure proves the  claim.
\end{proof}

If $\wis{C}$ is an S-monoid, it follows from Corollary~1 that every point in the topos $\widehat{\wis{C}}$ is isomorphic to some positive cone $\mathbb{Q}_+(s)$ where $s$ satisfies the condition of the above theorem. However, if $\wis{C}$ is only an s-monoid, other 'exotic' points exists and are described in Theorem~2.

\vskip 2mm

The previous theorem excludes the case of prime interest, that of $\mon$, but we easily deduce \cite[Thm. 2.2.(ii)]{CC}:

\begin{theorem}[Connes-Consani] The set of isomorphism classes of the arithmetic site $\widehat{\mon}$ is in natural bijection with the ad\`ele classes
\[
\mathbb{Q}^{\ast}_+ \backslash \mathbb{A}^f_{\mathbb{Q}} / \hat{\mathbb{Z}}^{\ast} \]
\end{theorem}

\begin{proof} For any $s \in \mathbb{S}$ and all natural numbers $n$ and $m$ dividing $s$,  $lcm(m,n)$ also divides $s$. This means that the rank one condition for $\mathbb{Q}_+(s)$ is always satisfied when the $S$-monoid is $\mon$. Therefore, isomorphism classes of points of $\widehat{\mon}$ are in natural bijection to the set of all equivalence classes $\mathbb{S}/\sim$.
\end{proof}

\section{A topology on $\mathbb{Q}^{\ast}_+ \backslash \mathbb{A}^f_{\mathbb{Q}} / \hat{\mathbb{Z}}^{\ast}=\mathbb{S}/\sim$}

We can turn the set of all supernatural numbers $\mathbb{S}$ into a compact Hausdorff topological space. Identify $\mathbb{S}$ with $\prod_{p \in \mathbb{P}} (\mathbb{N} \cup \{ \infty \})$ and take the product topology where we view $\mathbb{N} \cup \{ \infty \}$ as the one-point compactification of the discrete topology on $\mathbb{N}$, that is, a basis of open sets in $\mathbb{N} \cup \{ \infty \}$ is given by the singletons $\{ n \}$ for all $n \in \mathbb{N}$ and the sets $U_n = \{ m \geq n \} \cup \{ \infty \}$. This is the induced topology on $\mathbb{S}$ coming from the compact (resp. locally compact) topologies on $\widehat{\mathbb{Z}}$ and $\mathbb{A}^f_{\mathbb{Q}}$ via the identification $\mathbb{S} \leftrightarrow \widehat{\mathbb{Z}}/\widehat{\mathbb{Z}}^*$.

\vskip 2mm

However, the induced topology on the equivalence classes $\mathbb{S}/\sim$ is trivial, that is, the only non-empty open set is $\mathbb{S}/\sim$ itself. Indeed let $U = \prod U_p$ be an open set of $\mathbb{S}$, that is, each $U_p$ open in $\mathbb{N} \cup \{ \infty \}$, and suppose that $U$ is closed under the equivalence relation. 

If $s \in U$ with $s_p \in \mathbb{N}$, then all $t \in \mathbb{S}$ with $s_q=t_q$ for all primes $q \not= p$ and $s_q \in \mathbb{N}$ must also lie in $U$, whence $U_p = \mathbb{N} \cup \{ \infty \}$. If $s_p = \infty$ then $U_p$ being open implies there must be $t \in U$ with $t_p \in \mathbb{N}$ and we can repeat the foregoing, so also in this case $U_p$ must be $\mathbb{N} \cup \{ \infty \}$, so $U= \mathbb{S}$. 

\vskip 2mm

Also the SGA4-topology on $\widehat{\mon}$ is trivial, see \cite[IV.8.4.3]{SGA4}. Here, opens correspond to subobjects of the terminal object $\wis{1}$. Note that there is no terminal object in the small category $\mon$, but there is a terminal object in the presheaf topos $\widehat{\mon}$, given by the functor
\[
\wis{1}~:~\mon \rTo \wis{Sets} \qquad \bullet \mapsto \{ \ast \} \]
where $\{ \ast \}$ is a singleton with trivial right $\mon$-action. As the functor $S$ corresponding to the right $\mon$-set $\Omega_{\mon}$ is the sub-object classifier in $\widehat{\mon}$, see  \cite[I.4]{MM}, sub-objects of $\wis{1}$ correspond to natural transformations $1 \rTo S$, that is, to right $\mon$-maps
\[
\{ \ast \} \rTo \Omega_{\mon} \]
so the image must be a sieve which is an $\mon$-fixpoint. There are precisely two such sieves: $\emptyset$ and $\mon$.

\vskip 2mm

\begin{definition} The {\em sieve-topology} on $\wis{points}(\widehat{\mon})$ is the topology having as a basis of open sets the subsets
\[
\mathbb{X}(\wis{S}) = \wis{points}(\widehat{\mon}) \cap \wis{points}(\widehat{\wis{S}}) \]
where $\wis{S} \in \Omega_{\mon}$.
\end{definition}

This indeed defines a topology, for if $\wis{S} = c_1.\mon \cup c_2.\mon \cup \hdots$ and $\wis{S'} = d_1. \mon \cup d_2.\mon \cup \hdots$, we can define the s-monoid $\wis{S}.\wis{S'} = \cup_{i,j} lcm(c_i,d_j).\mon$. From Theorem~3 we then obtain that
\[
\mathbb{X}(\wis{S}) \cap \mathbb{X}(\wis{S'}) = \mathbb{X}(\wis{S}.\wis{S'}) \]
Hence, the sets $\mathbb{X}(S)$ form a basis of opens.

If we define for every $n \in \mon$ the open $\mathbb{X}(n) = \mathbb{X}(n.\mon) = \{ [\mathbb{Q}_+(s)]~:~n^{\infty} | s \}$, then we get for all $m,n \in \mon$
\[
\mathbb{X}(n) \cap \mathbb{X}(m) = \mathbb{X}(m.n) \quad \text{and} \quad \mathbb{X}(n) \cup \mathbb{X}(m) = \mathbb{X}(n.\mon \cup m.\mon) \]
Even though every sieve is of the form $\wis{S} = n_1.\mon \cup n_2.\mon \cup \hdots$, it is {\em not} true that the opens $\mathbb{X}(n)$ form a basis for the sieve-topology on the finite ad\`ele classes. In general we have for infinite unions
\[
\bigcup_{i =1}^{\infty} \mathbb{X}(n_i) \subsetneq \mathbb{X}(n_1 \mon \cup n_2 \mon \cup \hdots) \]
For example, take $n_i = p_i^{e_i}$ for distinct primes $p_i$, then with $s = \prod_i p_i^{e_i}$ we have that $[\mathbb{Q}_+(s)]$ belongs to the open set on the right hand side, but not to the union on the left hand side. That is, the topology requires an uncountable basis of open sets.

\begin{remark} If $p=[\mathbb{Q}_+(s)] \in \mathbb{X}({\wis{S}})$, and if $s|s'$, then also $q=[\mathbb{Q}_+(s')] \in \mathbb{X}(\wis{S})$, so we will not be able to separate points corresponding to positive cones whenever the corresponding supernatural numbers divide each other, or even, if they {\em weakly divide} each other. By this we mean
\[
s || s' \quad \text{if there exist $t,t' \in \mathbb{S}$ such that $s \sim t$, $s' \sim t'$ and $t | t'$} \]
In particular $[\mathbb{Q}_+(s_{\infty})]=[\mathbb{Q}_+]$ for $s_{\infty} = \prod_{p \in \mathbb{P}} p^{\infty}$ belongs to each of the open sets $\mathbb{X}(\wis{S})$.
\end{remark}

For this reason, we will say that the points $[\mathbb{Q}_+(s)]$ and $[\mathbb{Q}_+(s')]$ are {\em incomparable} if and only neither $s || s'$ nor $s' || s$. This is equivalent to saying that, with $\infty(t) = \{ p~:~p^{\infty} | t \}$,  either the sets $\infty(s)$ and $\infty(s')$ are incomparable, or that the subsets of primes $I = \{ p \in \mathbb{P}~:~\infty > s_p > s'_{p} \}$ and $J= \{ q \in \mathbb{P}~:~\infty > s'_{q} > s_q \}$ are both infinite. In the first case we can easily separate $s$ and $s'$, in the second case we take the sieves (or s-monoids)
\[
\wis{S} = \bigcup_{p \in I} p^{s_p}.\mon \quad \text{and} \quad \wis{S'} = \bigcup_{q \in J} q^{s'_q}.\mon \]
then we have that
\[
[\mathbb{Q}_+(s)] \in \mathbb{X}(\wis{S})-\mathbb{X}(\wis{S'}) \quad \text{and} \quad [\mathbb{Q}_+(s')] \in \mathbb{X}(\wis{S'})-\mathbb{X}(\wis{S}) \]
That is, the sieve-topology on the finite ad\`ele classes satisfies the $T_1$ separation property with respect to incomparable points.

\begin{theorem} The sieve-topology on the finite ad\`ele classes
\[
\mathbb{S}/\sim~\simeq \wis{points}(\widehat{\mon}) \simeq \mathbb{Q}^{\ast}_+ \backslash \mathbb{A}^f_{\mathbb{Q}} / \hat{\mathbb{Z}}^{\ast} \]
satisfies the following properties:
\begin{enumerate}
\item{It is compact.}
\item{It has an uncountable basis of opens.}
\item{Every non-empty open set is dense.}
\item{It satisfies the $T_1$ separation property for incomparable points.}
\end{enumerate}
\end{theorem}

\begin{proof}
For (1) recall that the only s-monoid $\wis{S}$ for which $[\mathbb{Q}_+(1)]=[\mon]$ is a point is $\mon$ itself. For (3) note that the point $[\mathbb{Q}_+(s_{\infty})]=[\mathbb{Q}_+]$ lies in each and every of the open sets $\mathbb{X}(\wis{S})$, so any two opens have a non-empty intersection. (2) and (4) were shown above.
\end{proof}

\vskip 5mm

\noindent
{\bf Acknowledgement : } I like to thank Pieter Belmans for explaining me that the sieve-topology is {\em not} the SGA4-topology on the topos $\widehat{\mon}$.

 \end{document}